\documentclass[11pt,oneside]{article}
\usepackage{amsfonts}
\usepackage{amsmath}
\usepackage{amsthm}
\usepackage[mathscr]{euscript}
\usepackage{mathrsfs}
\usepackage{indentfirst}
\usepackage{xcolor}

\newcommand{\lbl}[1]{\label{#1}}

\newtheorem{theo}{Theorem}[section]
\newtheorem{prop}{Proposition}[section]
\newtheorem{lem}{Lemma}[section]

\newcommand{\be}{\begin{equation}}
\newcommand{\ee}{\end{equation}}
\newcommand\bes{\begin{eqnarray}} \newcommand\ees{\end{eqnarray}}
\newcommand{\bess}{\begin{eqnarray*}}
\newcommand{\eess}{\end{eqnarray*}}

\newcommand\ep{\varepsilon}
\newcommand\kk{\left}
\newcommand\rr{\right}
\newcommand\dd{\displaystyle}

\newcommand\qq{\eqref}
\newcommand\yy{\infty}

\newcommand\mR{\mathbb{R}}
\newcommand\ol{\overline}
\newcommand\ud{\underline}

\begin{document}
\setlength{\baselineskip}{17pt} \pagestyle{myheadings}

\begin{center}{\Large\bf Dynamics for the diffusive Leslie-Gower model}\\[2mm]
{\Large\bf with double free boundaries}\footnote{This work was
supported by NSFC Grant 11371113}\\[4mm]
 {\Large  Mingxin Wang\footnote{Corresponding author. {\sl E-mail}: mxwang@hit.edu.cn}, \ \ Qianying Zhang}\\[0.5mm]
{\small Department of Mathematics, Harbin Institute of Technology, Harbin 150001, PR China}
\end{center}

\begin{quote}
\noindent{\bf Abstract.} In this paper we investigate a free boundary problem for the diffusive Leslie-Gower prey-predator model with double free boundaries in one space dimension. This system
models the expanding of an invasive or new predator species in which the free boundaries represent expanding fronts of the predator species. We first prove the existence, uniqueness and regularity of global solution. Then provide a spreading-vanishing dichotomy, namely the predator species either successfully spreads to infinity as $t\to\infty$ at both fronts and survives in the new
environment, or it spreads within a bounded area and dies out in the long run. The long time behavior of $(u,v)$ and criteria for spreading and vanishing are also obtained. Because the
term $v/u$ (which appears in the second equation) may be unbounded when $u$ nears zero,
it will bring some difficulties for our study.


 \end{quote}

 \section{Introduction}
 \setcounter{equation}{0} {\setlength\arraycolsep{2pt}

Prey-predator systems (or consumer-resource systems) are basic differential equation models for
describing the interactions between two species with a pair of positive-negative feedbacks.
The classical Leslie-Gower prey-predator model is (\cite{LG})
 \bess
 \left\{\begin{array}{lll}
 \dd\frac{{\rm d}u}{{\rm d}t}=u(a-u)-buv,\\[2.5mm]
 \dd\frac{{\rm d}v}{{\rm d}t}=\mu v\kk(1-\frac vu\rr),
 \end{array}\right.
 \eess
where $a$, $b$ and $\mu$ are positive constants, $u(t)$ and $v(t)$ represent
the population densities of prey and predator, respectively. In this model,
the prey is assumed to grow in logistic patterns. It is known that this system
has a globally asymptotically stable
equilibrium $\big(\frac{a}{1+b}, \frac{a}{1+b}\big)$.

The diffusive Leslie-Gower prey-predator model with homogeneous Neumann boundary conditions takes the form
 $$
 \left\{\begin{array}{lll}
 u_t-u_{xx}=u(a-u)-buv, \ \ &t>0, \ x\in\Omega,\\[1mm]
 v_t-dv_{xx}=\dd\mu v\kk(1-\frac vu\rr), \ \ &t>0, \ x\in\Omega,\\[1.5mm]
 \dd\frac{\partial u}{\partial\nu}=\frac{\partial u}{\partial\nu}=0, \ \ &t>0, \ x\in\partial\Omega,\\[1.5mm]
 u(0,x)=u_0(x)>0, \ \ v(0,x)=v_0(x)>0,\ \ &x\in\Omega,
 \end{array}\right.\eqno(P)
 $$
 where $\Omega$ is a bounded domain of $\mathbb R^N$. The global stability of $\big(\frac{a}{1+b}, \frac{a}{1+b}\big)$ for the problem (P) had been
 studied by many authors, see \cite{DHsu} for example.

In the problem (P), it is assumed that the habitats of prey and predator are the
same and fixed, and no flux through the boundary. However, in some situations,
predator and/or prey will have a tendency to emigrate from the
boundary to obtain their new habitat, i.e., they will move outward along
the unknown curve (free boundary) as time increases.
The spreading and vanishing of multiple species is an important content in understanding ecological complexity. In order to study the spreading and
vanishing phenomenon, many mathematical models have been established.

We assume that the prey distributes in the whole line $\mathbb{R}$ and the
predator exists initially in a bounded interval and invades into the new
environment from two sides. In such a situation the diffusive
Leslie-Gower prey-predator model with double
free boundaries can be written as
 \bes
 \left\{\begin{array}{lll}
 u_t-u_{xx}=u(a-u)-buv,\ \ &t>0,\ \ x\in\mR,\\[1.5mm]
 v_t-dv_{xx}=\dd\mu v\kk(1-\frac vu\rr),\ \ &t>0, \ \ g(t)<x<h(t),\\[1.5mm]
 v(t,x)\equiv 0, &t\ge 0,\ \ x\not\in(g(t), h(t)),\\[0.5mm]
 g'(t)=-\beta v_x(t,g(t)), \ &t\ge 0,\\[0.5mm]
 h'(t)=-\beta v_x(t,h(t)), \ &t\ge 0,\\[0.5mm]
 u(0,x)=u_0(x),\ \ &x\in\mR,\\[0.5mm]
v(0,x)=v_0(x),\ \ &-h_0\le x\le h_0, \\[0.5mm]
 g(0)=-h_0, \ \ h(0)=h_0,
 \end{array}\right.\lbl{1.1}
 \ees
where $a, b, d, h_0, \mu$ and $\beta$ are given positive constants. The initial functions $u_0(x),v_0(x)$ satisfy
 \bes
 u_0\in C_b(\mathbb{R}),\ u_0>0\ \ {\rm in} \ \mathbb{R};\ \ \
  v_0\in W^2_p((-h_0,h_0)),\ v_0(\pm h_0)=0,\ v_0>0\ \ {\rm in} \ (-h_0,h_0),
 \lbl{1.2}\ees
where $p>3$, $C_b(\mathbb{R})$ is the space of continuous and bounded functions
in $\mathbb{R}$. The free boundary condition $h'(t)=-\beta v_x(t,h(t))$ is the
Stefen type, and the deduction can refer to \cite{BDK} and \cite{WZ15}.

Free boundary problems of the classical Lotka-Volterra type prey-predator models had been investigated systematically by many authors, please refer to \cite{WZ15} (with double free boundaries), \cite{Wjde14, Wcnsns15, WZhang} (with homogeneous Dirichlet (Neumann, Robin) boundary conditions at the left side and free boundary at the right side), and \cite{ZWna14} (the prey distributes in the whole space $\mathbb{R}^N$, while the predator exists initially in a ball and invades into the new environment).

There were many related works for the classical Lotka-Volterra type competition models.
Authors of \cite{DL2, DWZ15, ZhaoW} investigated a competition model in which the invasive species exists initially in a ball and invades into the new environment, while the resident species distributes in the whole space $\mathbb{R}^N$. In \cite{GW12, WZjdde, ZhaoW}, two competition species are assumed to spread along the same free boundary at the right side and with homogeneous Dirichlet (Neumann, Robin) boundary conditions at the left side. Especially, the growth rates permit sign-changing in \cite{ZhaoW}. For the heterogeneous time-periodic environments, authors of \cite{CLW} and \cite{WZzamp} investigated the case with sign definite coefficients and the case with sign-changing growth rates, respectively.

The classical Lotka-Volterra type competition systems and prey-predator systems with different free boundaries had been studied in \cite{GW15, Wcom15, WZdd, Wu15}.

Without the predator in the environment (namely in the case $v\equiv 0$), \eqref{1.1}
reduces to a free boundary problem for $u$ considered in the pioneer work \cite{DLin}.
In this relatively simpler situation a spreading-vanishing dichotomy is known, and when
spreading happens, the spreading speed has been determined through a semi-wave problem
involving a single equation. More general results in this direction can be found in
\cite{CLZ, DGP, DLiang, DLou, KY, Wjde15, Wjfa16}, where \cite{CLZ} concerns with a nonlocal reaction term, \cite{DGP, Wjfa16} considers time-periodic environment,
\cite{DLiang} studies space-periodic environment, \cite{DLou, KY} investigates more general reaction terms. Particularly, in \cite{Wjde15, Wjfa16} the growth rates are allowed to change signs.

Free boundary problems of reaction diffusion equations and systems with advection
had been studied by many authors, refer to \cite{GKLZ, GLZ, KM, MWu, SLZ, WZZ, ZZL} for example.

This paper is organized as follows. In Section 2 we study the global existence, uniqueness, regularity and some estimates of $(u,v, g, h)$. Section 3 is concerned with the long time behaviors of $(u,v)$, and Section 4 deals with the criteria governing spreading and vanishing. Because the term $v/u$ may be unbounded when $u$ nears zero, it will bring some difficulties for the study.

At last we mention that for the free boundary problem of Holling-Tanner prey-predator model with double free boundaries
 \bess
 \left\{\begin{array}{lll}
 u_t-u_{xx}=u(a-u)-\dd\frac{buv}{m+u},\ \ &t>0,\ \ x\in\mR,\\[2.5mm]
 v_t-dv_{xx}=\dd\mu v\kk(1-\frac vu\rr),\ \ &t>0, \ \ g(t)<x<h(t),\\[2mm]
 v(t,x)\equiv 0, &t\ge 0,\ \ x\not\in(g(t), h(t)),\\[0.5mm]
 g'(t)=-\beta v_x(t,g(t)), \ &t\ge 0,\\[0.5mm]
 h'(t)=-\beta v_x(t,h(t)), \ &t\ge 0,\\[0.5mm]
 u(0,x)=u_0(x),\ \ &x\in\mR,\\[0.5mm]
v(0,x)=v_0(x),\ \ &-h_0\le x\le h_0, \\[0.5mm]
 g(0)=-h_0, \ \ h(0)=h_0,
 \end{array}\right.\lbl{1.1x}
 \eess
the methods used here are valid and the corresponding results are still true.

\section{Existence, uniqueness, regularities and estimates of global solution}
\setcounter{equation}{0} {\setlength\arraycolsep{2pt}

Set
  \[g^*=-\beta v_0'(-h_0),  \ \ h^*=-\beta v_0'(h_0), \ \ \mathbb{R}_+=(0,\infty),
 \ \  \ol{\mathbb{R}}_+=[0,\infty).\]
Then $g^*\le 0$ and $h^*\ge 0$. In order to facilitate the writing, we denote
  \[\Lambda=\big\{a, \, b, \, d, \, \mu, \, \beta, \, \alpha, \, p\}\]
with $0<\alpha<1-3/p$. For the given interval $I\subset\overline{\mathbb{R}}_+$, we set
 \bess
 I\times[g(t),h(t)]=\bigcup_{t\in I}\{t\}\times[g(t),h(t)], \ \
 I\times(g(t),h(t))=\bigcup_{t\in I}\{t\}\times(g(t),h(t)).\eess

\begin{theo}\lbl{th2.1}\, For any given $(u_0,v_0)$ satisfying $(\ref{1.2})$, the problem $(\ref{1.1})$ admits a unique global solution $(u, v, g, h)$ and
 \bess
 &u\in C_b(\overline{\mathbb{R}}_+\times\mR)\cap C^\yy(\mathbb{R}_+\times\mathbb{R}),\ \ 0<u\le\max\{a,\,\max u_0(x)\}:=A \ \ \mbox{in}\ \ \mathbb{R}_+\times\mathbb{R},&\\[.5mm]
 &v\in C^\yy(\mathbb{R}_+\times[g(t),h(t)]),\ \ 0<v\le\max\{A,\,\max v_0(x)\}:=B \ \ \mbox{in}\ \ \mathbb{R}_+\times(g(t),h(t)),&\\[.5mm]
 &g, h\in C^\yy(\mathbb{R}_+),\ \
 g'(t)<0, \ h'(t)>0 \ \ \mbox{in}\ \ \mathbb{R}_+.&
  \eess
Moreover, for any given $0<T<\yy$ and $0<\alpha<1-3/p$,
  $$v\in W_p^{1,2}((0,T)\times(g(t),h(t)))\hookrightarrow C^{\frac{1+\alpha}2,1+\alpha}([0,T]\times[g(t),h(t)]), \ \ g,\, h\in C^{1+\frac \alpha 2}([0,T]),$$
and
  \bes
 \|v\|_{W_p^{1,2}((0,T)\times(g(t),h(t))}+\|g,\,h\|_{C^{1+\frac \alpha 2}([0,T])}\leq C,\lbl{2.1}
  \ees
where the positive constant $C$ depends only on $T, \Lambda$, $h_0$, $g^*$, $h^*$, $\|u_0\|_{L^\yy(\mathbb R)}$, $ \|v_0\|_{W_p^2((-h_0,h_0))}$.
 \end{theo}

Before giving the proof of Theorem \ref{th2.1}, we first state a lemma which can be proved by the same way as that of \cite[Lemma 3.1]{WZ15} and the details will be omitted.

\begin{lem} $($Comparison principle$)$\label{lm2.1}\, Let $c$ and $T_0$ be two positive constants, $g_i, h_i\in C^1([0,T_0))$ and $g_i(t)<h_i(t)$ in $[0,T_0)$, $i=1,2$. Let $v_i\in C^{1,2}((0,T_0)\times(g_i(t),h_i(t)))$ and $v_{ix}\in C([0,T_0)\times[g_i(t),h_i(t)])$, $i=1,2$. Assume that $(v_1, g_1, h_1)$ and $(v_2, g_2, h_2)$ satisfy
$$
 \left\{\begin{array}{ll}
  v_{1t}-d v_{1xx}\leq\mu v_1(1-cv_1),\ \ &0<t<T_0,\ \ g_1(t)<x<h_1(t),\\[1mm]
 v_1(t,g_1(t))=0,\ \ g_1'(t)\geq-\beta v_{1x}(t,g_1(t)),&0<t<T_0,\\[1mm]
 v_1(t,h_1(t))=0,\ \ h_1'(t)\leq-\beta v_{1x}(t,h_1(t)),\ \ &0<t<T_0
 \end{array}\right.$$
and
 $$\left\{\begin{array}{ll}
  v_{2t}-d v_{2xx}\geq\mu v_2(1-cv_2),\ \ &0<t<T_0,\ \ g_2(t)<x<h_2(t),\\[1mm]
 v_2(t,g_2(t))=0,\ \ g_2'(t)\leq-\beta v_{2x}(t,g_2(t)),&0<t<T_0,\\[1mm]
 v_2(t,h_2(t))=0,\ \ h_2'(t)\geq-\beta v_{2x}(t,h_2(t)),\ \ &0<t<T_0,
 \end{array}\right.$$
respectively. If $g_1(0)\geq g_2(0), \, h_1(0)\leq h_2(0)$, $v_i(0,x)\geq 0$ on $[g_i(0),h_i(0)]$, and $v_1(0,x)\leq v_2(0,x)$ on $[g_1(0),h_1(0)]$,
then we have
  \[g_1(t)\geq g_2(t),\ h_1(t)\leq h_2(t) \ \ {\rm on}\  [0,T_0); \ \
 v_1(t,x)\leq v_2(t,x)\ \  {\rm on}\ \  [0,T_0)\times[g_1(t),h_1(t)].\]
  \end{lem}

{\bf Proof of Theorem \ref{th2.1}}. The proof will be divided into three steps.

{\it Step 1: Local existence and uniqueness}. The idea of this part comes from \cite{WZdd} and \cite{WZ15}. Let
 \bess
 y&=&\frac{2x-g(t)-h(t)}{h(t)-g(t)}, \\[.5mm]
 w(t,y)&=&u\big(t,\frac 12[(h(t)-g(t))y+h(t)+g(t)]\big), \\[.5mm]
 z(t,y)&=&v\big(t,\frac 12[(h(t)-g(t))y+h(t)+g(t)]\big).
 \eess
Then \eqref{1.1} is equivalent to
  \bes
 &&\left\{\begin{array}{lll}
 w_t-\rho(t)w_{yy}-\zeta(t,y)w_y=w(a-w-b z),\ &t>0, \ \ y\in\mR,\\[1mm]
 w(0,y)=u_0(h_0y):=w_0(y), \ &y\in\mR,\end{array}\right.
  \lbl{2.2}\\[2mm]
 &&\left\{\begin{array}{lll}
 z_t-d\rho(t)z_{yy}-\zeta(t,y)z_y=\dd\mu z\kk(1-\frac zw\rr),\ &t>0, \ \ -1<y<1,\\[1.5mm]
 z(t,\pm 1)=0,\ \ &t\ge 0, \\[0.5mm]
 z(0,y)=v_0(h_0y):=z_0(y), \ &-1\leq y\le 1,\end{array}\right.
  \lbl{2.3}\\[2mm]
 &&\left\{\begin{array}{lll}
 g'(t)=-\dd\beta\frac 2{h(t)-g(t)}z_y(t,-1),\ \ \ &t\ge 0,\\[3mm]
 h'(t)=-\dd\beta\frac 2{h(t)-g(t)}z_y(t,1),\ \ \ &t\ge 0,\\[3mm]
 g(0)=-h_0, \ \ h(0)=h_0,
 \end{array}\right.\lbl{2.4}
 \ees
where
 \[\rho(t)=\frac 4{(h(t)-g(t))^2}, \quad \zeta(t,y)=\frac{h'(t)+g'(t)}{h(t)-g(t)}+
 \frac{h'(t)-g'(t)}{h(t)-g(t)}y.\]
For $0<T\le\frac{h_0}{2(2+|g^*|+h^*)}:=T_1$,
we denote $I_T=[0,T]\times[-1,1]$, and
 $$ \begin{array}{l}
 \mathcal{D}^1_T=\big\{z\in C(I_T): \ z(0,y)=z_0(y), \ z(t,\pm 1)=0, \
  \|z-z_0\|_{C(I_T)}\leq 1\big\},\\[1mm]
\mathcal{D}^2_T=\{g\in C^1([0,T]):\, g(0)=-h_0, \ g'(0)=g^*, \
\|g'-g^*\|_{C([0,T])}\leq 1\},\\[1mm]
\mathcal{D}^3_T=\{h\in C^1([0,T]):\, h(0)=h_0, \ h'(0)=h^*, \
\|h'-h^*\|_{C([0,T])}\leq 1\}.
 \end{array} $$
Clearly, $\mathcal{D}_T= \mathcal{D}^1_T\times\mathcal{D}^2_T\times\mathcal{D}^3_T$
is a closed convex set of $C(I_T)\times[C^1([0,T])]^2$.

Next, we shall apply the contraction mapping theorem to show the existence and uniqueness result. Due to the choice of $T$, we see that, for
$(g,h)\in \mathcal{D}_T^2\times\mathcal{D}_T^3$,
 $$ |g(t)+h_0|+|h(t)-h_0|\leq T(\|g'\|_{C([0,T])}+\|h'\|_{C([0,T])})\leq {h_0}/2,$$
which implies
 \[h(t)-g(t)\ge h_0.\]
Given $(z, g, h)\in \mathcal{D}_T$, we set $z=0$ in $[0,T]\times((-\yy,-1]\cup[1,\yy))$ and substitute $(z(t,y), g(t), h(t))$ into \eqref{2.2}. Then \eqref{2.2} is a Cauchy problem of $w$. The standard theory (cf. \cite{LSU}) guarantees that the problem \eqref{2.2} admits a unique solution $w\in C([0,T]\times\mathbb{R})$. As $u_0(x)>0$ in $[-2h_0,2h_0]$, by use of the structure of $\mathcal{D}_T$ and the continuity of $w$, we can find a $T_2>0$ depending on $a$, $b$, $h_0$, $g^*$, $h^*$, $T_1$ and $u_0(x)$ such that
 \[\min_{[0,T]\times[-2,2]}w(t,y):=\delta(T)\ge\frac 12\min_{[-2,2]}w_0(y)=\frac 12\min_{[-2h_0,2h_0]}u_0(x)>0\]
provided $0<T\le T_2$.

Substituting this known function $w(t,y)$ into \eqref{2.3} and taking advantage of the $L^p$ theory and Sobolev's imbedding theorem we have that the problem \qq{2.3} admits a unique solution, denoted by $\tilde z(t, y)$, and $\tilde z\in W^{1,2}_p(I_T)\cap C^{\frac{1+\alpha}2,1+\alpha}(I_T)$,
 \bess
 \|\tilde z\|_{W^{1,2}_p(I_T)}\leq C_1, \ \ \|\tilde z\|_{C^{\frac{1+\alpha}2,1+\alpha}(I_T)}\leq C_1.
 \eess
for some constant $C_1>0$, depending on $d$, $\mu$, $p$, $\alpha$, $h_0$, $g^*$, $h^*$,
$\|v_0\|_{W_p^2((-h_0,h_0))}$ and $\frac 12\min_{[-2h_0,2h_0]}u_0(x)$.

Define
 \bess
  \tilde g(t)&=&-h_0-2\beta\int_0^t\frac 1{h(t)-g(t)}\tilde z_y(\tau,-1){\rm d}\tau,\\[1mm]
  \tilde h(t)&=&h_0-2\beta\int_0^t\frac 1{h(t)-g(t)}\tilde z_y(\tau,-1){\rm d}\tau.\eess
Then $\tilde g', \tilde h'\in C^{\frac \alpha 2}([0,T])$,
and
 \bess \|\tilde g'\|_{C^{\frac \alpha 2}([0,T])},\
 \|\tilde h'\|_{C^{\frac\alpha 2}([0,T])}\leq C_2, \eess
where $C_2$ depends on $C_1$, $\beta$, $h_0$, $g^*$ and $h^*$.

\vskip 4pt Now we define a mapping $\mathcal{F}: \mathcal{D}_T\rightarrow C(I_T)\times[C^1([0,T])]^2$ by
 \[\mathcal{F}(z, g, h)=(\tilde z,\tilde g,\tilde h).\]
Clearly, $(z,g, h)\in \mathcal{D}_T$ is a fixed point of $\mathcal{F}$ if and only if
$(w, z, g, h)$ solves (\ref{2.2})-\eqref{2.4}. Similar to the arguments in the proof of \cite[Theorem 2.1]{WZ15}, it can be shown that $\mathcal{F}$
maps $\mathcal{D}_T$ into itself and is a contraction mapping on $\mathcal{D}_T$ when $0<T\ll 1$, where $T$ depends only on $\Lambda$, $h_0$, $g^*$, $h^*$, $\|u_0\|_{L^\yy(\mathbb R)}$, $ \|v_0\|_{W_p^2((-h_0,h_0))}$ and $\frac 12\dd\min_{[-2h_0,2h_0]}u_0(x)$.
Thus $\mathcal{F}$ has a unique fixed point $(z,g, h)$ in $\mathcal{D}_T$ by the contraction mapping theorem, and $(z,g,h)\in C^{\frac{1+\alpha} 2,\,1+\alpha}(I_T)\times[C^{1+\frac\alpha 2}([0,T])]^2$. That is,  (\ref{2.2})-\eqref{2.4} have a unique solution $(w,z,g, h)$. By use of the $L^p$ and Schauder theories we can show that
 \[w\in C^{1+\frac\alpha2,2+\alpha}([\tau,T]\times[-L, L]), \ \
 z\in C^{1+\frac\alpha2,2+\alpha}([\tau,T]\times[-1, 1]),\ \
 g, h\in C^{1+\frac{1+\alpha}2}([\tau,T]),\]
for any given $0<\tau<T$ and $L>0$ (cf. \cite[Theorem 2.1]{Wjfa16}). Moreover,
$w,z>0$ by the maximum principle, and $z_y(t,-1)>0$, $z_y(t,-1)<0$ by the Hopf
boundary lemma, the latter imply that $g'(t)<0, h'(t)>0$ in $(0,T]$.
Therefore, the problem (\ref{1.1}) admits a unique solution $(u,v,g, h)$ and
 \[(u,v,g, h)\in C([0,T]\times\mathbb{R})\times C^{\frac{1+\alpha} 2,\,1+\alpha}([0,T]\times[g(t), h(t)])\times[C^{1+\frac\alpha 2}([0,T])]^2.\]

{\it Step 2: Global existence}. We extend the solution $(u,v,g, h)$ of \eqref{1.1} to the maximal time interval $[0, T_0)$ and show that $T_0=\yy$.

Firstly, by the maximum principle,
 \bess
 &0<u\le\max\{a,\,\max u_0(x)\}:=A\ \ \ \mbox{in} \ \ [0,T_0)\times\mathbb{R},&\\[.5mm]
 &0<v\le\max\{A,\,\max v_0(x)\}:=B\ \ \ \mbox{in} \ \ [0,T_0)\times(g(t),h(t)).&\eess
Using \qq{2.3}, \qq{2.4} and the Hopf boundary lemma we have that $g'(t)<0$ and $h'(t)>0$ in $[0,T_0)$.

Assume on the contrary that $T_0<\yy$. Let $(\bar v,\bar g, \bar h)$ be the unique solution of the following free boundary problem
\bess
 \left\{\begin{array}{lll}
 \bar v_t-d\bar v_{xx}=\mu \bar v,\ \ &t>0, \ \ \bar g(t)<x<\bar h(t),\\[0.5mm]
 \bar v(t,\bar g(t))=\bar v(t,\bar h(t))=0,\ \ \ &t\ge0,\\[0.5mm]
 \bar g'(t)=-\beta\bar v_x(t,\bar g(t)), \ &t\ge0,\\[0.5mm]
 \bar h'(t)=-\beta\bar v_x(t,\bar h(t)), \ &t\ge0,\\[0.5mm]
 \bar v(0,x)=v_0(x),\ \ &-h_0\le x\le h_0, \\[0.5mm]
 \bar g(0)=-h_0, \ \ \bar h(0)=h_0.
 \end{array}\right.
 \eess
In view of Lemma \ref{lm2.1},
 \[g(t)\ge\bar g(t)\ge\bar g(T_0), \ \ h(t)\le\bar h(t)\le\bar h(T_0) \ \ \ \mbox{in} \ \ [0,T_0).\]
Because $u$ satisfies
 \bess
 \left\{\begin{array}{lll}
 u_t-du_{xx}\ge u(a-bB-u),\ \ &0<t<T_0,\ \ x\in\mR,\\[0.5mm]
 u_x(t,0)=0,\ \ \ &0\le t<T_0,\\[0.5mm]
 u(0,x)=u_0(x)>0,\ \ &x\in\mR,
 \end{array}\right.
 \eess
we have
 \bes
 \min_{[0,T_0)\times[\bar g(T_0),\,\bar h(T_0)]}u(t,x)=\delta(T_0)>0.\lbl{2.5}\ees
Notice that $[g(t),h(t)]\subset[\bar g(T_0),\,\bar h(T_0)]$ for all $t\in[0,T_0)$,
we can find a positive constant $K(T_0)$ such that
 \[\mu v(1-v/u)\le K(T_0)\ \ \ \ \mbox{in} \, \ [0,T_0)\times[g(t),h(t)].\]
Similar to the proof of \cite[Lemma 2.1]{WZ15} we can show that there exists a positive constant $M(T_0)$, which depends only on $\Lambda$, $K(T_0)$, $\dd\min_{[0,h_0]}v_0'(x)$ and $\dd\max_{[-h_0,0]}v_0'(x)$, such that $g'(t)\ge-M(T_0)$, $h'(t)\le M(T_0)$ in $[0,T_0)$.\vskip 4pt

Let $0<\alpha<1-3/p$ and
  \[\Lambda(T_0)=\big\{\alpha,\, p,\, \bar g(T_0),\, \bar h(T_0),\, K(T_0),\, M(T_0)\big\}.\]
Applying the $L^p$ theory to \eqref{2.3} and the embedding theorem we have
 \[z\in W^{1,2}_p([0,T_0)\times(-1, 1))\hookrightarrow C^{\frac{1+\alpha}2,1+\alpha}([0,T_0]\times[-1, 1]),\]
and there exists a positive constant $C_1(T_0)$ depending only on $\Lambda(T_0)$ and $d$
such that
 \[\|z\|_{C^{\frac{1+\alpha} 2,\,1+\alpha}([0,T_0]\times[-1, 1])}\leq C_1(T_0).\]
Thus, by use of \qq{2.4}, $g, h\in C^{1+\frac\alpha 2}([0,T_0])$ and
 \[\|g, h\|_{C^{1+\frac\alpha 2}([0,T_0])}\leq C_2(T_0),\]
where $C_2(T_0)$ depends on $\Lambda(T_0)$, $C_1(T_0)$ and $\beta$.
Take advantage of the Schauder theory to \eqref{2.3} we have
 \[z\in C^{1+\frac\alpha2,2+\alpha}([T_0/2,T_0-\ep]\times[-1, 1]),
 \ \ \forall \ 0<\ep<T_0/2,\]
and there exists a positive constant $C_3(T_0)$, which depends only on $C_1(T_0)$ and $C_2(T_0)$, but not on $\ep$, such that
 \[\|z\|_{C^{1+\frac\alpha2,2+\alpha}([T_0/2,T_0-\ep]\times[-1, 1])}
 \leq C_3(T_0), \ \ \forall \ 0<\ep<T_0/2.\]
Then, in view of \qq{2.4}, we have $g, h\in C^{1+\frac{1+\alpha}2}([T_0/2,T_0-\ep])$ and
 \[\|g, h\|_{C^{1+\frac{1+\alpha}2}([T_0/2,T_0-\ep])}\leq C_4(T_0), \ \ \forall \ 0<\ep<T_0/2,\]
where $C_4(T_0)$ depends on $C_2(T_0)$ and $C_3(T_0)$.
Therefore, $v(t,\cdot)\in C^2([g(t),h(t)])$ for any
$T_0/2\le t<T_0$, and
 \[\|v(t,\cdot)\|_{C^2([g(t),h(t)])}+|g(t)|+|g'(t)|+|h(t)|+|h'(t)|
 \leq C_3(T_0), \ \ \forall \ \tau\le t<T_0.\]

Repeating the discussion of Step 1 we can find a positive constant $T$ depending only on
$a$, $b$, $d$, $\mu$, $\beta$, $\alpha$, $p$, $A$, $C_2(T_0)$ and $\delta(T_0)$ which
was given by \qq{2.5}, such that the solution of \qq{1.1} with initial
time $T_0-T/2$ can be extended uniquely to the time $T_0-T/2+T$.  But this
contradicts the definition of $T_0$.

{\it Step 3.} The regularity and estimate \qq{2.1} can be proved by the similar way to that of
\cite[Theorem 1.2]{Wzhaoy} and the details are omitted here. The proof is finished. \ \ \ \ \fbox{}

Let $(u,v,g,h)$ be the unique global solution of \qq{1.1}. As $g'(t)<0$ and $h'(t)>0$, we can define the limits $\dd\lim_{t\to\yy}g(t)=g_\yy\ge-\yy$ and $\dd\lim_{t\to\yy}h(t)=h_\yy\le\yy$.

\vskip 4pt At last, we shall give the uniform estimates of $v$ and $g',\,h'$ when $h_\yy-g_\yy<\yy$. To this purpose, we first state a proposition.

\begin{prop}\label{p2.1} \ Let $d, m, \theta, k, T$ be constants and $d, m,\,\theta>0$, $k,\,T\ge 0$. For any given $\ep, L>0$, there exist $T_\ep>T$ and $l_\ep>\max\big\{L,\frac{\pi}2\sqrt{d/m}\big\}$, such that when the function $w\in C^{1,2}((T,\yy)\times(-l_\ep,\,l_\ep))$ and satisfies $w\ge 0$,
 $$ \left\{\begin{array}{ll}
  w_t-d w_{xx}\geq\, (\leq)\, w(m-\theta w), \ \ &t>T, \ \ -l_\ep<x<l_\ep,\\[0.5mm]
 w(T,x)>0, &-l_\ep<x<l_\ep,
 \end{array}\right.$$
and for $t>T$, $w(t,\pm l_\ep)\geq\, (\leq)\, k$ if $k>0$, while $w(t,\pm l_\ep)\geq\, (=)\, 0$ if $k=0$, we must have
 \[ w(t,x)>m/\theta-\ep \ \ \big(w(t,x)<m/\theta+\ep\big), \ \ \ \forall \ t\geq T_\ep, \ \ x\in[-L,L].\]
This implies
 \[\liminf_{t\to\infty}w(t,x)\geq m/\theta-\ep \ \ \kk(\limsup_{t\to\infty}w(t,x)<m/\theta+\ep\rr)
  \ \ \ \mbox{uniformly\, on }\, [-L,L].\]
 \end{prop}

 This proposition can be proved by the same way as that of \cite[Proposition B.1]{WZ15} and the details are omitted here.

\begin{theo}\lbl{th2.2}\, Assume $b<1$. Let $(u,v,g,h)$ be the unique global solution of \qq{1.1}. Then, for any given $L>0$, there exists a positive constant $\sigma(L)>0$, which does not depend on $\beta$, such that $u(t,x)\ge\sigma(L)$ in $[0,\yy)\times[-L,L]$.
 \end{theo}

{\bf Proof}. It is easy to see from the first equation of \qq{1.1} that $\dd\limsup_{t\to\yy}\max_{x\in\mR}u(t,x)\le a$. For any given $\ep>0$, there exists $T_\ep\gg 1$ such that $u(t,x)\le a+\ep$ for all $t\ge T_\ep$ and $x\in\mR$. Then $v$ satisfies
\bess
 \left\{\begin{array}{lll}
 v_t-dv_{xx}\le\dd\mu v\kk(1-\frac v{a+\ep}\rr),\ \ &t\ge T_\ep, \ \ g(t)<x<h(t),\\[1.5mm]
 v(t,x)=0, &t\ge T_\ep, \ \ x\not\in(g(t), h(t)).
 \end{array}\right.
 \eess
Let $w(t)$ be the unique positive solution of
 \[w'(t)=\mu w\kk(1-\frac w{a+\ep}\rr), \ \ t>T_\ep; \quad w(T_\ep)
 =\max_{[g(T_\ep),h(T_\ep)]}v(T_\ep, x).\]
Then $\dd\lim_{t\to\yy}w(t)=a+\ep$, and $v(t,x)\le w(t)$ for all $t\ge T_\ep$, $x\in\mR$ by the comparison principle. Therefore,
 \[\limsup_{t\to\yy}\max_{x\in\mR}v(t,x)\le a\]
by the arbitrariness of $\ep$.

Take $0<\omega, \, \ep\ll 1$ such that $a-b(a+\omega)>0$. Then there exists $T>0$ such that
  \[v(t,x)\leq a+\omega,\ \ \forall \ t\geq T, \ x\in\mR.\]
Let $T_\ep$ and $l_\ep$ be given by Proposition \ref{p2.1} with $d=\theta=1$, $m=a-b(a+\omega)$,
$k=0$. It is clear that
 \bess
 \left\{\begin{array}{lll}
  u_t-u_{xx}\geq u(a-b(a+\omega)-u),\ &t\geq T,
   \ x\in[-l_\ep,l_\ep],\\[1mm]
   u(t,\pm l_\ep)>0,\ \ &t\geq T.
 \end{array}\right.
 \eess
Take advantage of Proposition \ref{p2.1} it arrives at $\dd\liminf_{t\to\infty}u(t,x)\geq a-b(a+\omega)-\ep$ uniformly on $[-L,L]$. The arbitrariness of $\ep$ and $\omega$ imply that
$\dd\liminf_{t\to\infty}u(t,x)\geq a(1-b)>0$ uniformly on $[-L,L]$. Notice that $u(t,x)>0$ on $[0,\yy)\times[-L,L]$, we can find a constant $\sigma(L)>0$ such that
  \[u(t,x)\ge\sigma(L) \ \ \ \mbox{in } \ \ [0,\yy)\times[-L,L].\]
The proof is finished. \ \ \ \fbox{}

\begin{theo}\lbl{th2.3}\, Assume $b<1$. Let $(u,v,g,h)$ be the unique global solution of \qq{1.1}. If $h_\yy-g_\yy<\yy$, then there exists a positive constant $C$, depends only on $\Lambda$, $g_\yy$ and $h_\yy$ such that
  \bes
  \|v(t,\cdot)\|_{C^1([g(t),\,h(t)])}\leq C, \ \ \forall \ t\ge 0; \ \ \ \ \|g', h'\|_{C^{{\alpha}/2}([0,\yy))}\leq C.
  \lbl{2.6}\ees
 \end{theo}

{\bf Proof}. The condition $h_\yy-g_\yy<\yy$ implies $g_\yy>-\yy$ and $h_\yy<\yy$. By Theorem \ref{th2.2}, there exists a positive constant $\sigma>0$, which does not depend on $\beta$, such that $u(t,x)\ge\sigma$ in $[0,\yy)\times[g_\yy,h_\yy]$.
Follow the proof of \cite[Theorem 2.1]{Wjfa16} step by step we can get the estimate \qq{2.6}. The details are omitted here. \ \ \ \ \fbox{}

\section{Long time behavior of $(u,v)$}
\setcounter{equation}{0} {\setlength\arraycolsep{2pt}

This section concerns with the limits of $(u(t,x),\,v(t,x))$ as $t\to\infty$.

{\it Case 1}: $h_\yy-g_\infty<\infty$. In this case we shall prove that
  \bes
 \lim_{t\to\infty}\max_{g(t)\leq x\leq h(t)}v(t,x)=0,\lbl{3.1}\ees
and
 \bes
 \lim_{t\to\infty}u(t,x)=a \ \ \mbox{uniformly on the compact subset of } \ \mR.
 \lbl{3.2}\ees
For this purpose, we first give a proposition.

\begin{prop}\label{l7.2}\, Let $d$, $C$, $\mu$ and $\eta_0$ be positive constants, $w\in W^{1,2}_p((0,T)\times(0,\eta(t)))$ for some $p>1$ and any $T>0$, and $w_x\in C([0,\infty)\times(0,\eta(t)])$, $\eta\in C^1([0,\infty))$. If $(w,\eta)$ satisfies
  \bess\left\{\begin{array}{lll}
 w_t-d w_{xx}\geq -C w, &t>0,\ \ 0<x<\eta(t),\\[1.5mm]
 w\ge 0,\ \ \ &t>0, \ \ x=0,\\[1.5mm]
 w=0,\ \ \eta'(t)\geq-\mu w_x, \ &t>0,\ \ x=\eta(t),\\[1.5mm]
 w(0,x)=w_0(x)\ge,\,\not\equiv 0,\ \ &x\in (0,\eta_0),\\[1.5mm]
 \eta(0)=\eta_0,
 \end{array}\right.\lbl{7z.25}\eess
and
   \bess
  &\dd\lim_{t\to\infty} \eta(t)=\eta_\infty<\infty, \ \ \ \lim_{t\to\infty} \eta'(t)=0,&
 \label{7.22}\\[1mm]
& \|w(t,\cdot)\|_{C^1([0,\,\eta(t)])}\leq M, \ \ \forall \ t>1&\label{7.23}
 \eess
for some constant $M>0$. Then
  \[\lim_{t\to\infty}\,\max_{0\leq x\leq \eta(t)}w(t,x)=0.\]
\end{prop}

\noindent{\it Proof}.\, Firstly, by the maximum principle we have $w(t,x)>0$ for $t>0$ and $0<x<\eta(t)$. Follow the proof of \cite[Theorem 2.2]{WZjdde} word by word we can prove this lemma and the details are omitted here. \ \ \ \fbox{}

\begin{theo}\lbl{th3.1}{\rm(}Vanishing{\rm)}\, Suppose $b<1$. Let $(u,v,g,h)$ be the unique global solution of \qq{1.1}. If $h_\yy-g_\infty<\infty$, then \qq{3.1} and \qq{3.2} hold.
\end{theo}

{\bf Proof}. Notice Theorem \ref{th2.3}, in view of Proposition \ref{l7.2} we can get $\dd\lim_{t\to\infty}\max_{0\leq x\leq h(t)}v(t,x)=0$ directly. Similarly, $\dd\lim_{t\to\infty}\max_{g(t)\leq x\leq 0}v(t,x)=0$. Thus \qq{3.1} holds. The proof of \qq{3.2} is standard and we shall omit the details. The interested readers can refer to the proof of \cite[Theorem 4.2]{WZ15}. \ \ \ \fbox{}

{\it Case 2}: $h_\yy-g_\infty=\infty$. In this case we shall prove that \bes
 \lim_{t\to\infty}u(t,x)=\frac{a}{1+b},\ \ \ \
\lim_{t\to\infty}v(t,x)=\frac{a}{1+b}
 \label{3.3}
 \ees
uniformly in any compact subset of $\mathbb{R}$. To do this we first show a proposition
which alleges that $h_\infty-g_\infty=\infty$ implies $g_\infty=-\infty$ and $h_\infty=\infty$.

\begin{prop}\label{p3.1} \ Assume $b<1$. Let $(u,v,g,h)$ be the unique global solution of \qq{1.1}. If $h_\infty-g_\infty=\infty$, then $g_\infty=-\infty$ and $h_\infty=\infty$.
\end{prop}

{\bf Proof.} Assume on the contrary that $h_\infty<\infty$. Then the condition  $h_\infty-g_\infty=\infty$ implies $g_\infty=-\infty$. There exists $T\gg 1$ such that
  \[h_0-g(T)>\pi\sqrt{d/\mu}.\]
Similar to the proof of Theorem \ref{th2.2}, we can find a constant $\sigma=\sigma(T)>0$ such that $u(t,x)\ge\sigma$ in $[T,\yy)\times[g(T), h_\yy]$. Then $v$ satisfies
  \bess
 \left\{\begin{array}{lll}
  v_t-dv_{xx}\dd\ge\mu v\kk(1-\frac v\sigma\rr),\ \ &t>T,\ \ g(T)<x<h(t),\\[1.5mm]
 v\ge 0, \ &t\ge T,\ \ x=g(T),\\[.5mm]
 v=0, \ h'(t)=-\beta v_x,\   \ &t\ge T,\ \ x=h(t).
 \end{array}\right.
 \eess
Let $(w,\eta)$ be the unique global solution of
  \bess
 \left\{\begin{array}{lll}
  w_t-dw_{xx}=\dd\mu w\kk(1-\frac w\sigma\rr),\ \ &t>T, \ \ g(T)<x<\eta(t),\\[0.5mm]
  w=0, \ &t\ge T,\ \ x=g(T), \\[0.5mm]
  w=0,\ \eta'(t)=-\beta w_x, \ &t\ge T,\ \ x=\eta(t),\\[0.5mm]
  w(T,x)=v(T,x),\ \ &g(T)\le x\le \eta(T), \\[0.5mm]
 \eta(T)=h(T).
 \end{array}\right.
 \eess
The comparison principle yields $h(t)\ge\eta(t)$ in $[T,\yy)$ and $v(t,x)\ge w(t,x)$
in $[T,\yy)\times[g(T), \eta(t)]$. Since
  \[\eta(T)-g(T)=h(T)-g(T)>h_0-g(T)>\pi\sqrt{d/\mu},\]
by use of \cite[Lemma 3.2]{Wjde15} we have that $\dd\lim_{t\to\yy}\eta(t)=\yy$, which implies $h_\yy=\yy$. A contradiction is obtained and so $h_\yy=\yy$. Similarly, we can show that $g_\yy=-\yy$.
The proof is complete.\ \ \ \ \fbox{}

\begin{theo}\label{th3.2}{\rm(}Spreading{\rm)}\, Assume $b<1$. Let $(u,v,g,h)$ be the unique global solution of \qq{1.1}. If
$g_\infty=-\infty$, $h_\infty=\infty$,
then \qq{3.3} holds.

\end{theo}

{\bf Proof.} This proof is similar to that of \cite[Theorem 4.3]{WZ15}. For the completeness we shall give the details.

Firstly,
   \bes\limsup_{t\to\infty}\max_{x\in\mR}u(t,x)\leq a,\ \ \limsup_{t\to\yy}\max_{x\in\mR}v(t,x)\le a, \label{3.4}\ees
see the proof of Theorem \ref{th2.2}.

Chosen $L\gg 1$ and $0<\omega, \, \ep\ll 1$ such that $a-b(a+\omega)>0$.
Let $l_\ep$ be given by Proposition \ref{p2.1} with $d=\theta=1$, $m=a-b(a+\omega)$,
$k=0$. Using (\ref{3.4}), we can choose a $T_1>0$ such that
 \[v(t,x)\leq a+\omega,\ \ \forall \ t\geq T_1, \ x\in[-l_\ep,l_\ep].\]
Then $u$ satisfies
 \bess
 \left\{\begin{array}{lll}
  u_t-u_{xx}\geq u[a-b(a+\omega)-u],\ &t\geq T_1,
   \ x\in[-l_\ep,l_\ep],\\[1mm]
   u(t,\pm l_\ep)>0,\ \ &t\geq T_1.
 \end{array}\right.
 \eess
In view of Proposition \ref{p2.1} we have $\dd\liminf_{t\to\infty}u(t,x)\geq a-b(a+\omega)-\ep$ uniformly on $[-L,L]$. In consideration of the arbitrariness of $L$, $\omega$ and $\ep$, it follows that
 \bes\liminf_{t\to\infty}u(t,x)\geq a-ba:=\ud a_1\ \
 \mbox{uniformly\, on\, the\, compact\,  subset\, of } \, \mathbb{R}. \label{3.6}\ees
By our assumption, $\ud a_1>0$.

Given $L\gg 1$ and $0<\omega, \ep\ll 1$. Let $l_\ep$ be determined by Proposition \ref{p2.1} with $m=\mu$, $\theta=\mu/(\ud a_1-\omega)$ and $k=0$. According to (\ref{3.6}) and $g_\infty=-\infty$, $h_\infty=\infty$, there is $T_2>0$ such that $u\geq\ud a_1-\omega$ in
$[T_2,\infty)\times[-l_\ep,l_\ep]$ and $g(t)<-l_\ep$, $h(t)>l_\ep$ for $t\geq T_2$. Consequently, $v$ satisfies
 \bess
 \left\{\begin{array}{lll}
  v_t-dv_{xx}\dd\geq\mu v\kk(1-\frac v{\ud a_1-\omega}\rr),\ &t\geq T_2,\ x\in[-l_\ep,l_\ep],\\[3mm]
  v(t,\pm l_\ep)\geq 0,\ \ &t\geq T_2.
 \end{array}\right.
 \eess
Similar to the above,
 \bes\liminf_{t\to\infty}v(t,x)\geq \ud{a}_1\ \ \mbox{uniformly\, on\, the\,
compact\,  subset\, of } \, \mathbb{R}. \label{3.7}\ees

Clearly, $a-b\ud a_1>0$. Take $L\gg 1$, $0<\omega, \ep\ll 1$ with $a-b(\ud a_1-\omega)>0$.
Let $l_\ep$ be given by Proposition \ref{p2.1} with
 \[d=\theta=1, \ \ m=a-b(\ud a_1-\omega), \ \ k=A:=\max\{a,\,\max u_0(x)\}.\]
By virtue of (\ref{3.7}) we can find a $T_3>0$ such that
 \[v(t,x)\geq \underline{a}_1-\omega,\ \ \
  \forall \ t\geq T_3, \ x\in[-l_\ep,l_\ep].\]
Thus $u$ satisfies
 \bess
 \left\{\begin{array}{lll}
  u_t-u_{xx}\leq u[a-b(\ud a_1-\omega)-u],
  \ &t\geq T_3, \ x\in[-l_\ep,l_\ep],\\[1mm]
   u(t,\pm l_\ep)\leq A,\ \ &t\geq T_3.
 \end{array}\right.
 \eess
The same as the above,
  \bes
  \limsup_{t\to\infty}u(t,x)\leq
a-b\ud a_1:=\bar{a}_1\ \ \mbox{uniformly\, on\, the\,
compact\,  subset\, of } \, \mathbb{R}.\lbl{xx.1}\ees

Given $L\gg 1$ and $0<\omega, \ep\ll 1$. Let $l_\ep$ be determined by Proposition \ref{p2.1} with
 \[m=\mu, \ \ \theta=\mu/(\bar a_1+\omega), \ \ k=B:=\max\{A,\,\max v_0(x)\},\]
where $A$ is given by the above. Thanks to \qq{xx.1}, there exists $T_4>0$ such that
 \[u(t,x)\leq \bar{a}_1+\omega,\ \ \
  \forall \ t\geq T_4, \ x\in[-l_\ep,l_\ep].\]
 Consequently, $v$ satisfies
 \bess
 \left\{\begin{array}{lll}
  v_t-dv_{xx}\dd\leq\mu v\kk(1-\frac v{\bar a_1+\omega}\rr),\ &t\geq T_4,\ x\in[-l_\ep,l_\ep],\\[3mm]
  v(t,\pm l_\ep)\leq B,\ \ &t\geq T_4.
 \end{array}\right.
 \eess
The same as the above,
  \[\limsup_{t\to\infty}v(t,x)\leq\bar{a}_1\ \ \mbox{uniformly\, on\, the\,
compact\,  subset\, of } \, \mathbb{R}.\]

Repeating the above procedure, we can find two sequences $\{\underline{a}_i\}_{i=1}^\yy$
and $\{\bar{a}_i\}_{i=1}^\yy$ such that, for all $i$,
 \bes
 \underline{a}_i\leq\liminf_{t\to\infty}u(t,x)\leq
 \limsup_{t\to\infty}u(t,x)\leq\bar{a}_i,\quad
 \underline{a}_i\leq\liminf_{t\to\infty}v(t,x)\leq
\limsup_{t\to\infty}v(t,x)\leq\bar{a}_i
\label{3.8}\ees
uniformly in the compact subset of $\mathbb{R}$. Moreover, these sequences can be
determined by the following iterative formulas:
 \[\ud{a}_1=a-ba, \ \ \bar a_i=a-b\ud a_i, \ \ \ud a_{i+1}=a-b\bar a_i,\ \ i=1,2,\cdots.\]
The direct calculation yields
 \[\bar a_1=a(1-b+b^2), \ \ \underline{a}_2=a(1-b+b^2-b^3), \ \  \bar{a}_2=a(1-b+b^2-b^3+b^4).\]
Using the inductive method we have the following expressions:
 \bess
 \ud a_i=a(1-b+b^2-\cdots+b^{2i-2}-b^{2i-1}),
 \ \  \bar{a}_i=a(1-b+b^2-\cdots+b^{2i}),
 \ \ \ i\geq 3.\eess
Because of $0<b<1$, one has
 \[\lim_{i\to\infty}\bar a_i=\lim_{i\to\infty}\underline{a}_i= a/(1+b).\]
This fact combined with (\ref{3.8}) allows us to derive \qq{3.3}.\ \ \ \ \fbox{}

\section{The criteria governing spreading and vanishing}
\setcounter{equation}{0}

We first give a necessary condition for vanishing.

\begin{theo}\label{th4.1}\, Assume $b<1$. Let $(u,v,g,h)$ be the unique global solution of \qq{1.1}. If
$h_\infty-g_\infty<\infty$,
then
  \bes
 h_\infty-g_\infty\leq\pi\sqrt{d/\mu}.\label{4.1}\ees
Hence, $h_0\geq \frac{\pi}2\sqrt{d/\mu}$ implies
$h_\infty-g_\infty=\infty$ due to $h'(t)-g'(t)>0$ for $t>0$.
\end{theo}

{\bf Proof.} The condition $h_\infty-g_\infty<\infty$ implies that
$\dd\lim_{t\to\infty}\|v(t,\cdot)\|_{C([g(t),h(t)])}=0$, $\dd\lim_{t\to\infty}u(t,x)=a$ uniformly in the compact subset of $\mathbb{R}$ (Theorem \ref{th3.1}). We assume $h_\infty-g_\infty>\pi\sqrt{d/\mu}$ to get a contradiction. For any given $0<\varepsilon\ll 1$, there exists $T\gg 1$ such that
  \bess
  &u(t,x)\dd\geq a-\varepsilon,\ \ \
  \forall \ t\geq T, \ x\in[g_\infty,h_\infty],&\\[1mm]
 &h(T)-g(T)>\pi\sqrt{d/\mu}.&\eess
Let $w$ be the unique solution of
  $$\left\{\begin{array}{ll}
 w_t-d w_{xx}=\dd\mu w\kk(1-\frac w{a-\varepsilon}\rr), \ \ &t>T, \ \, g(T)<x<h(T),\\[3mm]
 w(t,g(T))=w(t,h(T))=0, \ \ &t\ge T,\\[1mm]
  w(T,x)=v(T,x),&g(T)\le x\le h(T).
  \end{array}\right.$$
As $v$ satisfies
  $$\left\{\begin{array}{ll}
 v_t-d v_{xx}\ge\dd\mu v\kk(1-\frac v{a-\varepsilon}\rr), \ \ &t>T, \ \, g(T)<x<h(T),\\[3mm]
 v(t,g(T))\ge 0, \ \ v(t,h(T))\ge 0, \ \ &t\ge T,
  \end{array}\right.$$
the comparison principle gives $w\leq v$ in $[T,\infty)\times[g(T),\,h(T)]$. Since $h(T)-g(T)>
  \pi\sqrt{d/\mu}$,
it is well known that $w(t,x)\to \theta(x)$ as $t\to\infty$ uniformly on $[g(T),h(T)]$, where $\theta(x)$ is the unique positive solution of
  \[\left\{\begin{array}{ll}
 -d\theta_{xx}=\mu\theta\dd\kk(1-\frac \theta{a-\varepsilon}\rr)=0,\ \ &g(T)<x<h(T),\\[3mm]
 \theta(g(T))=\theta(h(T))=0.&\end{array}\right.\]
Hence, $\dd\liminf_{t\to\infty}v(t,x)\geq\lim_{t\to\infty}w(t,x)=\theta(x)>0$ in
$(g(T),h(T))$. This is a contradiction to (\ref{3.1}), and hence (\ref{4.1}) holds.
The proof is complete. \ \ \ \ \fbox{}

\vskip 2pt If $b<1$, by Theorem \ref{th4.1} and Proposition \ref{p3.1} we see that  $h_0\geq\frac{\pi}2\sqrt{d/\mu}$
implies $g_\infty=-\infty$ and $h_\infty=\infty$ for all $\beta>0$.

\vskip 2pt Now we discuss the case $h_0<\frac{\pi}2\sqrt{d/\mu}$.

\begin{lem}\label{lm4.1}\, Let $(u,v,g,h)$ be the unique global solution of \qq{1.1}. If $h_0<\frac{\pi}2\sqrt{d/\mu}$, then there exists $\beta_0>0$, depending on $d,h_0,\mu$ and $v_0(x)$, such that $g_\yy>-\yy$, $h_\infty<\infty$ provided $\beta\leq\beta_0$.
\end{lem}

{\bf Proof.} Obviously, $\lambda_1=\dd\frac{d}{4h_0^2}\pi^2$ and $\phi(x)=\dd\sin\frac{\pi(x+h_0)}{2h_0}$ are the principal eigenvalue and the corresponding positive eigenfunction of the following  problem
$$ \left\{\begin{array}{ll}
  -\phi_{xx}=\lambda\phi, \ \  -h_0<x<h_0,\\[1mm]
 \phi(\pm h_0)=0,
 \end{array}\right.$$
and there exists $k>0$ such that
 \[x\phi'(x)\le k\phi(x)\ \ \ \mbox{in} \ \ [-h_0, h_0].\]
The condition $h_0<\frac{\pi}2\sqrt{d/\mu}$ implies $\lambda_1>\mu$.
Let $0<\ep,\,\rho<1$ and $K>0$ be constants, which will be determined. Set
 \bess
 s(t)&=&1+2\ep-\ep {\rm e}^{-\rho t}, \ \ \eta(t)=h_0s(t), \ \ t\geq 0,\\[.5mm]
 w(t,x)&=&K{\rm e}^{-\rho t}\phi\left(x/s(t)\right), \ \ \ t\geq 0,\ \ -\eta(t)\leq x\leq \eta(t).
 \eess
Clearly, $w(t,\pm \eta(t))=0$.
Similar to the calculations in the proof of \cite[Lemma 3.4]{Wjde15} (\cite[Lemma 5.3]{Wjfa16}), we can show that there exists $\delta>0$ such that
 \[w_t-d w_{xx}-\mu w>0, \ \ \forall \ t>0, \ \ -\eta(t)<x<\eta(t)\]
for all $0<\ep, \rho\le\delta$ and all $K>0$. Fixed $0<\ep, \rho\le\delta$, then
 \[w(0,x)=K\phi\left(x/(1+\ep)\right)\ge v_0(x)\ \ \ \mbox{in} \ \ [-h_0, h_0]\]
provided $K\gg 1$. For these fixed $0<\ep, \rho\le\delta$ and $K\gg 1$, remember $\phi'(-h_0)>0$ and $\phi'(h_0)<0$, we can find a $\beta_0$: $0<\beta_0\ll 1$ such that, for all $0<\beta\le\beta_0$,
 \[-h_0\ep\rho\le-\beta\frac 1{s(t)}K\phi'(-h_0), \ \
 h_0\ep\rho\ge-\beta\frac 1{s(t)}K\phi'(h_0), \ \ \forall \ t\ge 0.\]
This implies
 \[-\eta'(t)\le-\beta w_x(t,-\eta(t)), \ \ \eta'(t)\ge-\beta w_x(t,\eta(t)), \ \ \forall \ t\ge 0.\]
Because of $v$ satisfies
 \[v_t-d v_{xx}-\mu v<0, \ \ \forall \ t>0, \ \ g(t)<x<h(t),\]
by the comparison principle we conclude
 \[g(t)\ge-\eta(t), \ \ h(t)\le\eta(t), \ \ \forall \ t\ge 0,\]
and so $g_\yy\ge -\eta(\yy)=-(1+2\ep)h_0$, $h_\yy\le \eta(\yy)=(1+2\ep)h_0$. The proof is complete.
\ \ \ \fbox{}

\begin{lem}\label{lm4.2}\, Let $C$ be a positive constant. For any given positive constants $h_0, L$, and any function $\bar v_0\in W^2_p((-h_0,h_0))$ with $p>1$, $\bar v_0(\pm h_0)=0$ and $\bar v_0>0$ in $(-h_0,h_0)$, there exists $\beta^0>0$ such that when $\beta\geq\beta^0$ and $(\bar v, \bar g, \bar h)$ satisfies
 \bes
 \left\{\begin{array}{ll}
   \bar v_t-\bar v_{xx}\geq -C \bar v, \ &t>0, \ \bar g(t)<x< \bar h(t),\\[0.5mm]
  \bar v=0, \ \bar g'(t)=-\beta \bar v_x, \ &t\geq 0, \ x=\bar g(t),\\[0.5mm]
 \bar v=0, \ \bar h'(t)=-\beta \bar v_x, \ &t\geq 0, \ x=\bar h(t),\\[0.5mm]
 \bar v(0,x)=\bar v_0(x),\ \ \ &-h_0\leq x\leq h_0,\\[.5mm]
 \bar g(0)=-h_0, \ \bar h(0)=h_0,
  \end{array}\right.\label{3.5}
 \ees
we must have $\dd\lim_{t\to\infty}\bar g(t)<-L$, $\dd\lim_{t\to\infty}\bar h(t)>L$.
\end{lem}

{\bf Proof}. Follow the proof of \cite[Lemma 3.2]{WZhang} step by step and use the
comparison principle, we can prove the conclusion.
The details are omitted here. \ \ \ \fbox{}

\begin{lem}\label{lm4.3} \, Assume that $b<1$. Let $(u,v,g,h)$ be the unique global
solution of \qq{1.1}. If $h_0<\frac{\pi}2\sqrt{d/\mu}$, then exists $\beta^0>0$ such
that $h_\yy-g_\yy>\pi\sqrt{d/\mu}$ for all $\beta\geq\beta^0$.
\end{lem}

{\bf Proof.} Write $(u^\beta, v^\beta, g^\beta, h^\beta)$ in place of $(u,v, g, h)$
to clarify the dependence of the solution of \qq{1.1} on $\beta$. Assume on the
contrary that there exist $\{\beta_n\}$ with $\beta_n\to\yy$ such that $h^{\beta_n}_\yy-g^{\beta_n}_\yy\le\pi\sqrt{d/\mu}$ for all $n$.
In view of Theorem \ref{th2.2}, there exists a positive constant $\sigma>0$
such that $u^{\beta_n}(t,x)\ge\sigma$ in $[0,\yy)\times[-\pi\sqrt{d/\mu},\,
\pi\sqrt{d/\mu}]$ for all $n$.
Therefore $\mu\big(1-\frac {v^{\beta_n}}{u^{\beta_n}}\big)\ge -C$ in $[0,\yy)
\times[-\pi\sqrt{d/\mu},\,\pi\sqrt{d/\mu}]$ for some positive constant $C$ and
all $n$, and so
 \[v^{\beta_n}_t-dv^{\beta_n}_{xx}\ge -Cv^{\beta_n},\ \ t>0, \ \ g^{\beta_n}(t)<x<h^{\beta_n}(t).\]
According to Lemma \ref{lm4.2}, there exists $\beta^0>0$ such that
$h^\beta_\yy-g^\beta_\yy>\pi\sqrt{d/\mu}$ for all $\beta\ge\beta^0$. Since $\beta_n\to\yy$, we
get a contradiction and the proof is complete. \ \ \ \fbox{}

Making use of Lemmas \ref{lm4.1} and \ref{lm4.3}, we can prove the following theorem by
the same manner as that of \cite[Theorem 5.2]{WZ15} and the details will be omitted.

\begin{theo}\label{th4.2} \, Assume that $b<1$. Let $(u,v,g,h)$ be the unique global solution of \qq{1.1}. If $h_0<\frac{\pi}2\sqrt{d/\mu}$, then there exist $\beta^*\geq\beta_*>0$, depending on $a$, $b$, $d$, $\mu$, $u_0(x)$, $v_0(x)$ and $h_0$, such that $g_\infty=-\infty$ and $h_\infty=\infty$ when $\beta>\beta^*$, and $h_\infty-g_\infty\leq\pi\sqrt{d/\mu}$ when  $\beta\leq\beta_*$ or $\beta=\beta^*$.
\end{theo}

Now we state the criteria for spreading ($g_\infty=-\infty$, $h_\infty=\infty$)
and vanishing ($h_\infty-g_\infty<\yy$).

\begin{theo}\lbl{th5.1}\, Assume that $b<1$. Let $(u,v,g,h)$ be the unique global solution of \qq{1.1}.

{\rm(i)}\, If $h_0\geq\frac{\pi}2\sqrt{d/\mu}$, then $g_\infty=-\infty$ and $h_\infty=\infty$ for all $\beta>0$;

{\rm(ii)}\, If $h_0<\frac{\pi}2\sqrt{d/\mu}$, then there exist $\beta^*\geq\beta_*>0$, such that $g_\infty=-\infty$ and $h_\infty=\infty$
for $\beta>\beta^*$, while $h_\infty-g_\infty\leq\pi\sqrt{d/\mu}$ for
$\beta\leq\beta_*$  or $\beta=\beta^*$.
\end{theo}

From the above discussion we immediately obtain the following spreading-vanishing dichotomy and criteria for spreading and vanishing.

\begin{theo}\lbl{th5.3} Assume $b<1$. Let $(u,v,g,h)$ be the unique global solution of \qq{1.1}. Then the following alternative holds:

Either\vspace{-2mm}\begin{itemize}
\item[{\rm(i)}]Spreading: $g_\yy=-\yy$, $h_\infty=\infty$ and
 \[\lim_{t\to\infty}u(t,x)=\frac{a}{1+b},\ \ \ \
\lim_{t\to\infty}v(t,x)=\frac{a}{1+b}\]
uniformly in any compact subset of $\mathbb{R}$,\vspace{-2mm}\end{itemize}
or
\vspace{-2mm}\begin{itemize}
\item[{\rm(ii)}]Vanishing: $h_\infty-g_\yy\leq\pi\sqrt{d/\mu}$ and
 \bess
 &\dd\lim_{t\to\infty}\max_{g(t)\leq x\leq h(t)}v(t,x)=0,&\\[1mm]
 &\dd\lim_{t\to\infty}u(t,x)=a \ \ \mbox{uniformly on the compact subset of} \ \mR.&
 \eess
\end{itemize}\vspace{-2mm}

Moreover,

{\rm(iii)} If $h_0\ge\frac 12\pi\sqrt{d/\mu}$, then $g_\yy=-\yy$, $h_\infty=\infty$ for all $\beta>0$;\vskip 4pt

{\rm(iv)}\, If $h_0<\frac 12\pi\sqrt{d/\mu}$, then there exist $\beta^*\geq\beta_*>0$, such that $g_\infty=-\infty$ and $h_\infty=\infty$ for $\beta>\beta^*$, while $h_\infty-g_\infty\leq\pi\sqrt{d/\mu}$ for $\beta\leq\beta_*$ or $\beta=\beta^*$.
\end{theo}

\section{Discussion}
\setcounter{equation}{0} {\setlength\arraycolsep{2pt}

In this paper we have examined the Leslie-Gower prey-predator model with double
free boundaries $x=g(t)$ and $x=h(t)$ for the predator. We envision that the
prey distributes in the whole space, while the predator initially occupy a finite
region $[-h_0, h_0]$ and invades into the new environment. The dynamics of
(\ref{1.1}) exhibits a spreading-vanishing dichotomy:

(i) When spreading happens, both prey and predator will stabilize at the unique
positive equilibrium state $\big(\frac{a}{1+b}, \frac{a}{1+b}\big)$ as $t\to\yy$.
This behavior is the same as that of the initial-boundary problem (P).

(ii) When vanishing occurs, the predator will spread within a bounded area and
dies out in the long run, the prey will stabilize at a positive equilibrium state.

The criteria governing spreading and vanishing indicate that both spreading and
vanishing are completely determined by the initial habitat of the predator and
initial densities of the prey and predator, and the moving parameter/coefficient
$\beta$ of free boundaries.

These results tell us that in order to control the prey species (pest species)
we should put predator species (natural enemies) at the initial state at least
in one of three ways: (i) expand the initial habitat of predator, (ii) increase
the moving parameter/coefficient of free boundaries, (iii) augment the initial
density of the predator species.

These theoretical results may be helpful in the prediction and prevention of
biological invasions.

\end{document}